\documentclass[a4paper,11pt]{amsart}

\usepackage{graphicx}
\usepackage{mathptmx}
\usepackage{amsmath}
\usepackage{amssymb}
\usepackage{enumitem}
\usepackage{xcolor}

\newmuskip\pFqmuskip

\newcommand*\pFq[6][8]{%
  \begingroup 
  \pFqmuskip=#1mu\relax
  \mathcode`=\string"8000
  \begingroup\lccode`\~=`\,
  \lowercase{\endgroup\let~}\pFqcomma
  F^{#2}_{#3}{\left(\genfrac..{0pt}{}{#4}{#5}\bigg|#6\right)}%
  \endgroup
}
\newcommand{\pFqcomma}{\mskip\pFqmuskip}

\newtheorem{theorem}{Theorem}[section]

\newtheorem{remark}[theorem]{Remark}

\begin{document}

\title[Moments of a random variable arising Laplacian random variable]{Moments of a random variable arising from Laplacian random variable}

\author{Taekyun  Kim}
\address{Department of Mathematics, Kwangwoon University, Seoul 139-701, Republic of Korea}
\email{tkkim@kw.ac.kr}
\author{Dae San  Kim}
\address{Department of Mathematics, Sogang University, Seoul 121-742, Republic of Korea}
\email{dskim@sogang.ac.kr}


\begin{abstract}
Let $X$ be the Laplacian random variable with parameters $(a,b)=(0,1)$, and let $(X_{j})_{j \ge 1}$ be a sequence of mutually independent copies of $X$. In this note, we explicitly determine the moments of the random variable $\sum_{k=1}^{\infty}\frac{X_{k}}{2k \pi}$ in terms of the Bernoulli and Euler numbers.
\end{abstract}

\maketitle

\section{Introduction} 
The Bernoulli numbers $B_{n}$ and the Euler numbers $E_{n}$ are respectively defined by 
\begin{equation}
\frac{t}{e^{t}-1}=\sum_{n=0}^{\infty}B_{n}\frac{t^{n}}{n!},\quad \frac{2}{e^{t}+1}=\sum_{n=0}^{\infty}E_{n}\frac{t^{n}}{n!},\quad (\mathrm{see}\ [1-4]). \label{1}
\end{equation}
The first few terms of $B_n$ are given by:
\begin{align}
&B_0=1,\,B_1=-\frac{1}{2},\,B_2=\frac{1}{6},\,B_4=-\frac{1}{30},\,B_6=\frac{1}{42},\,B_8=-\frac{1}{30},\,B_{10}=\frac{5}{66},\, \label{1-1} \\
&B_{12}=-\frac{691}{2730},\, B_{14}=\frac{7}{6},\,B_{16}=-\frac{3617}{510},\,B_{18}=\frac{43867}{798},\,B_{20}=-\frac{174611}{330},\dots; \nonumber \\ 
&B_{2k+1}=0,\,\,(k \ge 1).\nonumber
\end{align}
The first few terms of $E_n$ are given by:
\begin{align}
&E_0=1,\,E_1=-\frac{1}{2},\,E_3=\frac{1}{4},\,E_5=-\frac{1}{2},\,E_7=\frac{17}{8},\,E_9=-\frac{31}{2},\,E_{11}=\frac{691}{4},\, \label{1-2}\\
&E_{13}=-\frac{5461}{2},\, E_{15}=\frac{929569}{16},\,E_{17}=-\frac{3202291}{2},\,E_{19}=\frac{221930581}{4},\dots; \nonumber\\
& E_{2k}=0,\,\,(k \ge 1). \nonumber
\end{align} \par
A random variable $X$ is the Laplacian random variable with parameters $a$ and $b(>0)$, which is denoted by $X\sim L(a,b)$, if its probability density function is given by
\begin{equation}
f(x)=\frac{1}{2b}e^{-\frac{|x-a|}{b}},\ x\in(-\infty,\infty),\quad (\mathrm{see}\ [5,7]),\label{2}
\end{equation}
where $a$ is the local parameter and $b(>0)$ is the scale parameter. \par
The Euler's product expansion for the sine function is the identity 
\begin{equation}
\frac{\sin\pi x}{\pi x}=\prod_{j=1}^{\infty}\bigg(1-\frac{x^{2}}{j^{2}}\bigg),\quad \big(x\in (-\infty,\infty)\big),\quad (\mathrm{see}\ [2]). \label{3}
\end{equation}
This identity was used by Euler in 1735 to give a solution of the Basel problem. \par 
Let $X \sim L(0,1)$, and let $(X_{j})_{j \ge 1}$ be a sequence of mutually independent copies of $X$. In this note, we determine the moments of the random variable $Y=\sum_{k=1}^{\infty}\frac{X_{k}}{2k \pi}$. Indeed, we show that $E[Y^{2n}]=(-1)^{n}\big(\frac{2n}{2^{2n}}E_{2n-1}+B_{2n}\big)$,\,\, $(n \in \mathbb{N})$, and that all odd moments of $Y$ vanish (see Theorem 2.1).

\section{Moments of a random variable arising from Laplacian random variable} 
For $X\sim L(0,1)$, let us assume that $(X_{j})_{j\ge 1}$ is a sequence of mutually independent copies of the random variable $X$. From \eqref{2}, we note that 
\begin{align}
E\Big[e^{\frac{X_{k}}{2\pi k}t}\Big]&=\frac{1}{2}\int_{-\infty}^{\infty}e^{\frac{x}{2\pi k}t}e^{-|x|}dx \label{4}\\
&=\frac{1}{2}\bigg[\int_{-\infty}^{0}e^{x\big(1+\frac{t}{2\pi k}\big)}dx+\int_{0}^{\infty}e^{-\big(1-\frac{t}{2\pi k}\big)x}dx \nonumber \\
&=\frac{1}{2}\bigg[\frac{1}{1+\frac{1}{2\pi k}t}+\frac{1}{1-\frac{t}{2\pi k}	}\bigg]=\frac{1}{1-\big(\frac{t}{2\pi k}\big)^{2}}, \nonumber
\end{align}
where $k$ is a positive integer and $-2 \pi<t<2 \pi$. \par 
Thus, by \eqref{4}, we get 
\begin{align}
\prod_{k=1}^{\infty}\bigg(\frac{1}{1-\big(\frac{t}{2\pi k}\big)^{2}}\bigg)&=\prod_{k=1}^{\infty}E\Big[e^{\frac{X_{k}}{2\pi k}t}\Big]=E\bigg[\prod_{k=1}^{\infty}e^{\frac{X_{k}}{2\pi k}t}\bigg] \label{5} \\
&=E\bigg[e^{\sum_{k=1}^{\infty}\frac{X_{k}}{2\pi k}t}\bigg]. \nonumber
\end{align}
On the other hand, by \eqref{3}, we get 
\begin{align}
&\prod_{k=1}^{\infty}\bigg(\frac{1}{1-\big(\frac{t}{2\pi k}\big)^{2}}\bigg)=\cfrac{\frac{t}{2}}{\sin\frac{t}{2}}=\cfrac{\frac{t}{2}}{\frac{e^{\frac{it}{2}}-e^{-\frac{it}{2}}}{2i}}=i\frac{t}{e^{\frac{it}{2}}-e^{-\frac{it}{2}}}\label{6}\\
&=it\bigg(\frac{e^{\frac{it}{2}}-1+1}{e^{it}-1}\bigg)=\frac{it}{2}\bigg(\frac{2}{e^{\frac{it}{2}}+1}\bigg)+\frac{it}{e^{it}-1}\nonumber \\
&=\frac{it}{2}\sum_{n=0}^{\infty}E_{n}\cfrac{\big(\frac{it}{2}\big)^{n}}{n!}+\sum_{n=0}^{\infty}B_{n}\frac{(it)^{n}}{n!}\nonumber \\
&=\frac{it}{2}+\frac{it}{2}\sum_{n=1}^{\infty}E_{2n-1}\cfrac{\big(\frac{it}{2}\big)^{2n-1}}{(2n-1)!}+1-\frac{it}{2}+\sum_{n=1}^{\infty}\frac{B_{2n}}{(2n)!}(it)^{2n}\nonumber \\
&=1+\sum_{n=1}^{\infty}\frac{(-1)^{n}}{(2n-1)!}E_{2n-1}\bigg(\frac{t}{2}\bigg)^{2n}+\sum_{n=1}^{\infty}\frac{(-1)^{n}}{(2n)!}B_{2n}t^{2n} \nonumber \\
&=1+\sum_{n=1}^{\infty}(-1)^{n}\bigg(\frac{2n}{2^{2n}}E_{2n-1}+B_{2n}\bigg)\frac{t^{2n}}{(2n)!}. \nonumber
\end{align}
By \eqref{5} and \eqref{6}, we get 
\begin{align}
1+\sum_{n=1}^{\infty}(-1)^{n}\bigg(\frac{2n}{4^{n}}E_{2n-1}+B_{2n}\bigg)\frac{t^{2n}}{(2n)!}&=\prod_{k=1}^{\infty}\bigg(\cfrac{1}{1-\big(\frac{t}{2k\pi}\big)^{2}}\bigg)=E\Big[e^{\sum_{k=1}^{\infty}\frac{X_{k}}{2k\pi}t}\Big]\label{7} \\
&=	\sum_{n=0}^{\infty}E\bigg[\bigg(\sum_{k=1}^{\infty}\frac{X_{k}}{2k\pi}\bigg)^{n}\bigg]\frac{t^{n}}{n!}.\nonumber
\end{align}
Therefore, by comparing the coefficients on both sides of \eqref{7}, we obtain the following theorem. 
\begin{theorem}
For $X\sim L(0,1)$, let $(X_{j})_{j\ge 1}$ be a sequence of mutually independent copies of the random variable $X$. Then we have 
\begin{displaymath}
E\bigg[\bigg(\sum_{k=1}^{\infty}\frac{X_{k}}{2k\pi}\bigg)^{2n}\bigg]=(-1)^{n}\bigg(\frac{2n}{2^{2n}}E_{2n-1}+B_{2n}\bigg), 
\end{displaymath}
and 
\begin{displaymath}
E\bigg[\bigg(\sum_{k=1}^{\infty}\frac{X_{k}}{2k\pi}\bigg)^{2n-1}\bigg]=0,\quad (n\in\mathbb{N}). 
\end{displaymath}
\end{theorem}
\begin{remark} 
As is known, the Bernoulli and Euler numbers are related by:
\begin{equation}
E_{n}=-\frac{2(2^{n+1}-1)}{n+1}B_{n+1},\quad (n \ge 0). \label{8}
\end{equation}
For example, this follows from the equation (14) of [6].
Thus, from Theorem 2.1 and \eqref{8}, we have the following alternative expression:
\begin{align}
E\bigg[\bigg(\sum_{k=1}^{\infty}\frac{X_{k}}{2k\pi}\bigg)^{2n}\bigg]&=(-1)^{n-1}\Big(1-\frac{1}{2^{2n-1}}\Big)B_{2n} \label{9}\\
&=\Big(1-\frac{1}{2^{2n-1}}\Big)|B_{2n}|,\quad (n \in \mathbb{N}). \nonumber
\end{align}
Thus we have
\begin{equation*}
E\bigg[\bigg(\sum_{k=1}^{\infty}\frac{X_{k}}{2k\pi}\bigg)^{2n}\bigg] \sim |B_{2n}|,\,\, \textnormal{as}\,\,\, n \rightarrow \infty,\quad E\bigg[\bigg(\sum_{k=1}^{\infty}\frac{X_{k}}{2k\pi}\bigg)^{2n+1}\bigg]=|B_{2n+1}|,\quad (n\in\mathbb{N}). 
\end{equation*}
\end{remark}
Finally, we illustrate Theorem 2.1 by using \eqref{9}.
\begin{align*}
&E\bigg[\bigg(\sum_{k=1}^{\infty}\frac{X_{k}}{2k\pi}\bigg)^{2}\bigg] =\frac{1}{12},\quad
E\bigg[\bigg(\sum_{k=1}^{\infty}\frac{X_{k}}{2k\pi}\bigg)^{4}\bigg]=\frac{7}{240},\\
&E\bigg[\bigg(\sum_{k=1}^{\infty}\frac{X_{k}}{2k\pi}\bigg)^{6}\bigg]=\frac{31}{1344}, \quad
E\bigg[\bigg(\sum_{k=1}^{\infty}\frac{X_{k}}{2k\pi}\bigg)^{8}\bigg]=\frac{127}{3840},\\
&E\bigg[\bigg(\sum_{k=1}^{\infty}\frac{X_{k}}{2k\pi}\bigg)^{10}\bigg]=\frac{2555}{33792}, \quad
E\bigg[\bigg(\sum_{k=1}^{\infty}\frac{X_{k}}{2k\pi}\bigg)^{12}\bigg]=\frac{1414477}{5591040},\\
&E\bigg[\bigg(\sum_{k=1}^{\infty}\frac{X_{k}}{2k\pi}\bigg)^{14}\bigg]=\frac{57337}{49152},
\quad E\bigg[\bigg(\sum_{k=1}^{\infty}\frac{X_{k}}{2k\pi}\bigg)^{16}\bigg]=\frac{118518239}{16711680},\\
&E\bigg[\bigg(\sum_{k=1}^{\infty}\frac{X_{k}}{2k\pi}\bigg)^{18}\bigg]=\frac{5749691557}{104595456}, \quad
E\bigg[\bigg(\sum_{k=1}^{\infty}\frac{X_{k}}{2k\pi}\bigg)^{20}\bigg]
=\frac{91546277357}{173015040}.
\end{align*}

\end{document}